\documentclass[12pt]{article}
\usepackage{epsfig}
\usepackage{amsfonts}
\usepackage{amssymb}
\usepackage{latexsym}
\usepackage{amsmath,latexsym, amscd, amsthm}

\newtheorem{lemma}{Lemma}

\newtheorem{theorem}{Theorem}
\newtheorem{corollary}{Corollary}

\newtheorem{conjecture}{Conjecture}

\newenvironment{remark}[1][Remark]
           {\medbreak\noindent {\em #1. \enspace}}
           {\par \medbreak}

\makeatletter \@addtoreset{equation}{section} \makeatother

\def\ppt{\frac{\partial}{\partial t}}

\begin{document}

\title{Compact Gradient Shrinking Ricci Solitons \\with Positive Curvature Operator}
\author{Xiaodong Cao\\
\textit{cao@math.columbia.edu}}

\maketitle
\begin{abstract}
In this paper, we first derive several identities on a compact
shrinking Ricci soliton. We then show that a compact gradient
shrinking soliton must be Einstein, if it admits a Riemannian
metric with positive curvature operator and satisfies an integral
inequality. Furthermore, such a soliton must be of constant
curvature.
\end{abstract}

\section{Introduction and Main Theorems}
Hamilton started the study of the Ricci flow in \cite{H3}. In
\cite{HPCO}, Hamilton has classified all compact manifolds with
positive curvature operator in dimension four. Since then, the
Ricci flow has become a powerful tool for the study of Riemannian
manifolds, especially for those manifolds with positive curvature.
Perelman made significant progress in his recent work
~\cite{perelman1} and ~\cite{perelman2}.

Suppose we have a solution to the Ricci flow
\begin{align}
\ppt g_{ij}= -2R_{ij}
\end{align}
on a compact Riemannian manifold $M$ with Riemannian metric
$g(t)$. Ricci soliton emerges as the limit of the solutions of the
Ricci flow. A solution to the Ricci flow is called a Ricci soliton
if it moves only by a one-parameter group of diffeomorphism and
scaling. If the vector field which induce the diffeomorphism is in
fact the gradient of a function, we call it a gradient Ricci
soliton. For a gradient shrinking Ricci soliton, we have the
equation
\begin{align}
R_{ij}+\nabla_i \nabla_j f = \frac{1}{2\tau} g_{ij}~,
\end{align}
where $\tau = T-t$. $T$ is the time the soliton becomes a point,
and $f$ is called Ricci potential function. In the special case
when $f$ is a constant, then we have an Einstein manifold.

Besides the above equation, a gradient shrinking Ricci soliton
must also satisfies the following equations,
\begin{align}
R+ \triangle f=\frac{n}{2\tau}
\end{align}
and
\begin{align}
\label{eq:1.4} R+ |\nabla f|^2=\frac{f-c}{\tau}~,
\end{align}
where $c$ is a constant in space. The last equation (\ref{eq:1.4})
determines the value of $f$. The Ricci potential function $f$
satisfies the following evolution equation,
\begin{align}
\ppt f= |\nabla f|^2~.
\end{align}

Inspired by his own work in \cite{HPCO} and \cite{Hsurvey},
Hamilton made the following conjecture:

\begin{conjecture} (\textbf{Hamilton})
A compact gradient shrinking Ricci soliton with positive curvature
operator must be Einstein.
\end{conjecture}

On the other hand, it is a well-known theorem of Tachibana
\cite{tachibana} that any compact Einstein manifold with positive
sectional curvature must be of constant curvature. Hence
Hamilton's conjecture is a generalization of the Tachibana
theorem, since Einstein manifolds are special Ricci solitons with
constant Ricci potential functions.

In this paper, we first derive a sequence of identities on
gradient shrinking Ricci solitons. Then we show that the above
conjecture is in fact true provided that the Ricci soliton
satisfies an integral inequality.

One of our main theorems is the following:
\begin{theorem}
\label{theorem:1} Let $(M, g(t))$ be a compact gradient shrinking
Ricci soliton, then $M$ must be of constant curvature if its
curvature operator is positive and satisfies the following
inequality,
\begin{equation}
\label{eq:1.6} \frac{1}{2} \int |Rc|^2|\nabla f|^2 e^{-f}\leq \int
Ke^{-f}+\int R_{ijkl}R_{ik}f_j f_l e^{-f}~,
\end{equation}where
\begin{equation}
\label{eq:1.7} K=(\nabla_i \nabla_j R_{ik}-\nabla_j \nabla_i
R_{ik})R_{jk}~.
\end{equation}
\end{theorem}

In Section Two, we first derive some integral identities about
Riemannian curvature on gradient shrinking Ricci solitons. More
precisely, we prove the following two identities,

\begin{theorem} \label{theorem:2} On a compact gradient shrinking Ricci soliton, we
have
\begin{equation}
\label{eq:1.8} \int Rm(Rc,Rc)e^{-f}=\frac{1}{2\tau} \int |Rc|^2
e^{-f}+\frac{1}{2}\int |div~Rm|^2 e^{-f}
\end{equation} and
\begin{equation} \label{eq:1.9}
\int Rm(Rc,Rc)e^{-f}=\frac{1}{2\tau} \int |Rc|^2 e^{-f}+\int
|\nabla Rc|^2 e^{-f}-\frac{1}{2}\int |div~Rm|^2 e^{-f}~,
\end{equation}
where
\begin{equation} \label{eq:1.10}
Rm(Rc, Rc)=R_{ijkl}R_{ik}R_{jl}~.
\end{equation}
\end{theorem}
As a corollary of Theorem~\ref{theorem:2}, we have
\begin{corollary} \label{corollary:1} On a compact gradient shrinking Ricci soliton, we
have
\begin{equation} \label{eq:1.11}
\int |\nabla Rc|^2 e^{-f} =\int |div~ Rm|^2 e^{-f}~.
\end{equation}
Moreover, (\ref{eq:1.8}) and (\ref{eq:1.9}) can be written as
follows,
\begin{equation}
\label{eq:1.12} \int Rm(Rc,Rc)e^{-f}=\frac{1}{2\tau} \int |Rc|^2
e^{-f}+\frac{1}{2}\int |\nabla Rc|^2 e^{-f}~.
\end{equation}
\end{corollary}

In Section Three, we derive some identities about Ricci curvature,
i.e., we show the following theorem,
\begin{theorem} \label{theorem:3} On a compact gradient shrinking Ricci soliton, we
have
\begin{equation}
\label{eq:1.13} \frac{1}{2} \int |Rc|^2 \triangle (e^{-f})
=\frac{1}{2} \int |\nabla Rc|^2e^{-f} +\int Ke^{-f} +\int
R_{kljp}R_{kj}f_l f_p e^{-f}~.
\end{equation}
\end{theorem}

In Section Four, we prove Theorem 1 under the hypothesis of
positive curvature operator and inequality (\ref{eq:1.6}).

{\bf Acknowledgement:}We would like to thank Professor Gang Tian,
for first bringing this problem to our attention, and for many
valuable suggestions during this work. We would also like to thank
Professor Richard Hamilton, for his patience and guidance during
this work.


We are indebt to Professor Bennett Chow, who shared his own notes
in this direction with us. We would like to thank him for his
generous comments. We would also like to thank Professor Tom
Ilmanen and Professor Duong H. Phong for many helpful discussion.

\section{Identities of Riemannian Curvature}

In this section, we prove Theorem~\ref{theorem:2}. On a gradient
shrinking Ricci soliton, we have the following identities:
\begin{align}\label{eq:2.1} \nonumber
&(div~Rm)_{jkl} = R_{ijkl,i}=\nabla_i R_{ijkl}=\nabla_i R_{klij}
\\ \nonumber
 =& -\nabla_k R_{ijli}-\nabla_l R_{ijik}
= \nabla_k R_{jl}-\nabla_l R_{jk} \\ \nonumber =& -\nabla_k
f_{jl}+\nabla_l f_{jk}
= \nabla_{l}\nabla_{k}f_{j}-\nabla_{k}\nabla_{l}f_{j} \\
=& R_{lkjp}f_{p}~.
\end{align}
Hence we have the following two identities,
\begin{equation}
\nabla_i (R_{ijkl}e^{-f})=0
\end{equation}
and
\begin{equation}
\nabla_i (R_{ik}e^{-f})=0~.
\end{equation}
Using integration by parts, we derive that
\begin{align*}
&\int {|div~ R_m|^2 e^{-f}}\\
 =& \int{R_{lkjp}f_{p}(-R_{jk,l}+R_{jl,k}) e^{-f}} \\
=&\int {R_{lkjp}f_p R_{jl,k} e^{-f}} - \int {R_{lkjp}f_p  R_{jk,l}
e^{-f}}\\
=& -\int {R_{lkjp}f_{pk} R_{jl} e^{-f}} + \int{R_{lkjp}f_{pl}
R_{jk} e^{-f}} \\
=& -\int {R_{lkjp} R_{lj} f_{kp} e^{-f}} - \int{R_{kljp} R_{kj}
f_{lp} e^{-f}} \\
=& -2\int {R_{lkjp} R_{lj} f_{kp} e^{-f}}~.
\end{align*}
Hence we have the following lemma:
\begin{lemma} On a gradient shrinking Ricci soliton, we have
\begin{equation}
\label{eq:2.4} \int {R_{lkjp} R_{lj} f_{kp}
e^{-f}}=-\frac{1}{2}\int |div~Rm|^2e^{-f} \leq 0~.
\end{equation}
\end{lemma}
Now we can prove (\ref{eq:1.8}) in Theorem~\ref{theorem:2}.
\begin{proof}
By the above lemma and the gradient shrinking Ricci soliton
equation:
$$f_{kp}=\frac{1}{2\tau} g_{kp}-R_{kp}~,$$
we can derive
\begin{align*}
&\int {|div~ R_m|^2 e^{-f}}\\
=&-2\int {R_{lkjp} R_{lj} (\frac{1}{2\tau} g_{kp}-R_{kp})
e^{-f}}\\
=&-\frac{1}{\tau} \int |Rc|^2 e^{-f}+2\int Rm(Rc,Rc) e^{-f}~,
\end{align*}
so we have
$$\int Rm(Rc,Rc)e^{-f}=\frac{1}{2\tau} \int |Rc|^2
e^{-f}+\frac{1}{2}\int |div~Rm|^2 e^{-f}~.$$
\end{proof}
Before we prove (\ref{eq:1.9}), we first prove the following two
lemmas:
\begin{lemma}
\begin{equation}
\nabla_i \nabla_j R_{ik}-\nabla_j \nabla_i
R_{ik}=R_{jm}R_{mk}-R_{ijmk}R_{im}
\end{equation}
\end{lemma}
\begin{proof}
Using the formula
$$\nabla_i \nabla_j R_{lk}-\nabla_j \nabla_i
R_{lk}=-R_{ijml}R_{mk}-R_{ijmk}R_{lm}~,$$  and let $i=l$ in the
above formula and take the sum.
\end{proof}

\begin{lemma}On a gradient shrinking Ricci soliton,
\begin{equation}
-2\int \nabla_k R_{jl} \nabla_l R_{jk} e^{-f}=\frac{1}{\tau}\int
|Rc|^2 e^{-f}-2\int Rm(Rc,Rc)e^{-f}~.
\end{equation}
\end{lemma}
\begin{proof}
\begin{align*}
&-2\int \nabla_k R_{jl} \nabla_l R_{jk} e^{-f}\\
=&2\int R_{jk}(\nabla_i \nabla_j R_{ik}-\nabla_j R_{ik}
f_i)e^{-f}\\
=&2\int R_{jk}(\nabla_i \nabla_j R_{ik})e^{-f}-2\int
R_{jk}\nabla_j R_{ik}f_i e^{-f}\\
=&2\int R_{jk}(\nabla_j \nabla_i
R_{ik}+R_{mj}R_{mk}-R_{ijmk}R_{im}) e^{-f}+2\int
R_{ik}R_{jk}f_{ij}e^{-f}\\
=&0+2\int R_{jk}R_{mj}R_{mk}e^{-f}+ 2 \int
R_{ik}R_{jk}f_{ij}e^{-f}-2\int R_{ijmk}R_{im}R_{jk}e^{-f}\\
=&2\int R_{jk}R_{ki}(f_{ij}+R_{ij})e^{-f}-2\int R_{ijmk}R_{im}R_{jk}e^{-f}\\
=&\frac{1}{\tau} \int |Rc|^2 e^{-f}-2\int Rm(Rc,Rc)e^{-f}~.
\end{align*} This finishes the proof of the lemma.
\end{proof}
We used the following lemma in the above,
\begin{lemma}
\label{lemma:4} On a gradient shrinking Ricci soliton, we have
\begin{equation}
\int \nabla_j \nabla_i R_{ik}R_{jk}e^{-f}=0~.
\end{equation}
\end{lemma}
\begin{proof}
\begin{align*}
\int \nabla_j \nabla_i R_{ik}R_{jk}e^{-f} =-\int \nabla_i R_{ik}
\nabla_j( R_{jk}e^{-f})=0~.
\end{align*}
\end{proof}

Now we can prove (\ref{eq:1.9}) in Theorem~\ref{theorem:2}.
\begin{proof}
\begin{align*}
&\int {|div~ R_m|^2 e^{-f}}\\
=&\int {|\nabla_k R_{jl}-\nabla_l R_{jk}|^2 e^{-f}}\\
=&2\int |\nabla Rc|^2 e^{-f}-2\int {\nabla_k R_{jl} \nabla_l
R_{jk} e^{-f}}\\
=&2\int |\nabla Rc|^2 e^{-f}+\frac{1}{\tau} \int |Rc|^2
e^{-f}-2\int Rm(Rc,Rc)e^{-f}~,
\end{align*} so
$$\int Rm(Rc,Rc)e^{-f}=\frac{1}{2\tau} \int |Rc|^2
e^{-f}+\int |\nabla Rc|^2 e^{-f}-\frac{1}{2} \int {|div~ R_m|^2
e^{-f}}~.$$
\end{proof}

By (\ref{eq:1.8}) and (\ref{eq:1.9}) we have the
corollary~\ref{corollary:1}.

\section{Identities of Ricci Curvature}

Because of the soliton equation, there will be several identities
for Ricci curvature on the gradient shrinking Ricci solitons. We
first prove Theorem~\ref{theorem:3}. By using (\ref{eq:2.1}), we
derive that
\begin{equation}
\Delta R_{jk}=\nabla_i(\nabla_i R_{jk})=\nabla_i
(\nabla_jR_{ik}-R_{ijkl}f_l)=\nabla_i\nabla_jR_{ik}-(\nabla_iR_{ijkl})f_l-R_{ijkl}f_{li}~,
\end{equation}
so
\begin{equation} <\Delta Rc,
Rc>=\nabla_i\nabla_iR_{jk}R_{jk}=\nabla_i\nabla_jR_{ik} R_{jk}
-(\nabla_iR_{ijkl})f_l R_{jk}-R_{ijkl}f_{li}R_{jk}~,
\end{equation}
and
\begin{equation} \frac{1}{2} \Delta |Rc|^2=\frac{1}{2}\Delta
(R_{jk}R_{jk})=\nabla_i(\nabla_i R_{jk} R_{jk})=(\Delta R_{jk}
R_{jk})+|\nabla Rc|^2~.
\end{equation}

Furthermore, we have
\begin{align*}
 \frac{1}{2} \int \Delta |Rc|^2 e^{-f} =\frac{1}{2}\int
    |Rc|^2 \Delta e^{-f}
\end{align*}
so
\begin{align}\label{eq:3.4} \nonumber
&\frac{1}{2} \int |Rc|^2 \Delta e^{-f}\\ \nonumber =&\int < \Delta
Rc, Rc>e^{-f} + \int |\nabla Rc|^2e^{-f} \\ \nonumber
  =& \int |\nabla Rc|^2e^{-f} +\int (\nabla_i \nabla_jR_{ik}R_{jk}-\nabla_j \nabla_iR_{ik}R_{jk})e^{-f}
  \\ \nonumber
  &+ \int \nabla_j \nabla_i R_{ik}R_{jk}e^{-f} -\int \nabla_i R_{ijkl} f_l R_{jk}e^{-f} -\int
     R_{ijkl} f_{li} R_{jk}e^{-f}\\ \nonumber
  =& \int |\nabla Rc|^2e^{-f} +\int Ke^{-f}+\int \nabla_j \nabla_i R_{ik}R_{jk}e^{-f}-\int \nabla_i R_{ijkl} f_l R_{jk}e^{-f}
   -\int R_{ijkl} f_{li} R_{jk}e^{-f}\\
=& \int |\nabla Rc|^2e^{-f} +\int Ke^{-f}-\int \nabla_i R_{ijkl}
f_l R_{jk}e^{-f}-\int R_{ijkl} f_{li} R_{jk}e^{-f}~.
\end{align}
We used Lemma \ref{lemma:4} in the last equation.

Plug (\ref{eq:2.1}) and (\ref{eq:2.4}) into (\ref{eq:3.4}), apply
Corollary~\ref{corollary:1}, we obtain
\begin{align}\nonumber
\frac{1}{2} \int |Rc|^2 \Delta e^{-f}=&\int |\nabla Rc|^2e^{-f}
+\int Ke^{-f}-\int \nabla_i R_{ijkl} f_l R_{jk}e^{-f}
-\frac{1}{2}\int |\nabla Rc|^2e^{-f}\\
=& \frac{1}{2}\int |\nabla Rc|^2e^{-f} +\int Ke^{-f}+\int
R_{kljp}R_{jk}f_l f_p e^{-f}~.
\end{align}

If we assume that the metric on the gradient shrinking Ricci
soliton has positive curvature, then $$\int R_{kljp}R_{jk}f_p f_l
e^{-f}$$ is a positive term. In fact, this is true for any metric
with positive curvature operator. We have
\begin{lemma} Let $(M,g)$ be a Riemannian manifold with
positive curvature operator, then
$$R_{ijkl}R_{ik}f_jf_l \geq 0$$ point-wise.
\end{lemma}

\begin{proof}
$$R_{ikjl}=\sum_{\alpha} \lambda_{\alpha} \omega^{\alpha}_{ik}
\omega^{\alpha}_{jl}~,$$ where
$$\omega^{\alpha}=\omega^{\alpha}_{ik} dx^i \wedge dx^k$$ are
$2$-forms (in fact, they are the eigenfunctions of the curvature
operator). And $$\lambda_{\alpha}\geq 0~,$$ so
\begin{align*}
 R_{ikjl}R_{ij}f_k f_l=\sum_{\alpha} \lambda_{\alpha}
[\omega^{\alpha}_{ik} \omega^{\alpha}_{jl} R_{ij}f_k f_l]~,
\end{align*}
with
\begin{align*}
\label{eq:3.10} \omega^{\alpha}_{ik} \omega^{\alpha}_{jl}
R_{ij}f_k f_l=R_{ij} (\omega^{\alpha}_{ik}
f_k)(\omega^{\alpha}_{jl}f_l)=R_{ij}\gamma^{\alpha}_i
\gamma^{\alpha}_j \geq 0
\end{align*}
and $$\gamma^{\alpha}_i=\omega^{\alpha}_{ik}f_k~.$$
\end{proof}

\begin{remark} It's an easy calculation to see that
\begin{equation}\nonumber
\label{eq:3.11} \frac{1}{2} \int |Rc|^2 \Delta
e^{-f}=\frac{1}{\tau} \int Rc(\nabla f, \nabla f)e^{-f} \geq 0~.
\end{equation}
\end{remark}

\section{Proof of Theorem 1}
For a compact Riemannian manifold with positive curvature
operator, we first need the following lemma of Berger:
\begin{lemma} (\textbf{Berger})
Assume $T$ is a symmetric two tensor on a Riemannian manifold
$(M,g)$ with non-negative sectional curvature, then
$$K=(\nabla_i \nabla_j T_{ik}-\nabla_j \nabla_i
T_{ik})T_{jk} \geq 0~.$$ In fact, $$K=\sum_{i<j}
R_{ijij}(\lambda_i-\lambda_j)^2~,$$ where $\lambda_i$'s are the
eigenvalues of $T$.
\end{lemma}
We apply this lemma in the special case of $$T=Rc~.$$ Then we know
that our $K$ which is defined in (\ref{eq:1.7}) is non-negative.

By combining Theorem \ref{theorem:3} and inequality
(\ref{eq:1.6}), we can prove Theorem 1.
\begin{proof}
By
\begin{align*}
&\frac{1}{2}\int | Rc|^2 \triangle(e^{-f})\\
=&\int \nabla_i R_{jk} R_{jk}f_ie^{-f}\\
\leq& \frac{1}{2} \int |\nabla Rc|^2e^{-f}+\frac{1}{2}\int
|Rc|^2|\nabla f|^2e^{-f}\\
\leq& \frac{1}{2} \int |\nabla Rc|^2e^{-f}+\int Ke^{-f}+\int
R_{ijkl}R_{ik}f_j f_l e^{-f}~,
\end{align*}
we show that for all $i$, $j$ and $k$ we have
$$\nabla_i R_{jk}=R_{jk}f_i~,$$ and
$$ \int R_{kljp}R_{kj}f_l f_p e^{-f}=\int R_{kljp} \nabla_l R_{kj}
f_p e^{-f}=-\int R_{kljp} R_{kj}f_{lp}e^{-f}= \frac{1}{2}\int
|\nabla~Rc|^2 e^{-f}~.$$So $$\int Ke^{-f}= 0~,$$ hence $$K\equiv
0$$ and $f$ is a constant. Therefore, the soliton must be of
constant curvature.
\end{proof}

\nocite{*}
\bibliographystyle{plain}
\bibliography{bio}

\end{document}